\documentclass{article}[12pt]
\usepackage[english]{babel}
\usepackage[letterpaper,top=3cm,bottom=4cm,left=3cm,right=3cm,marginparwidth=1.75cm]{geometry}
\usepackage{amsmath}
\usepackage{amsfonts}
\usepackage{graphicx}
\usepackage[colorlinks=true, allcolors=blue]{hyperref}

\newtheorem{theorem}{Theorem}
\newtheorem{lemma}{Lemma}
\newtheorem{prop}{Proposition}

\title{The Eisenlohr-Farris algorithm for fully transitive polyhedra.}
\author{Eric Paulí Pérez Contreras}

\begin{document}
\maketitle

\begin{abstract}
The purpose of this note is to present a method for classifying three-dimensional polyhedra in terms of their symmetry groups. This method is constructive and it is described in terms of the conjugation classes of crystallographic groups in $\mathbb{E}^3$. For each class of groups $\Gamma$ the method can generate without duplication all polyhedra in three-dimensional space on which $\Gamma$ acts fully-transitively. It was proposed by J. M. Eisenlohr and S. L. Farris for generating every fully transitive polyhedra in $\mathbb{E}^d$. We also illustrate how the method can be applied in the euclidean space $\mathbb{E}^3$ by generating a new fully transitive polyhedron.
\end{abstract}

\section{Introduction}

\begin{figure}
\centering
\includegraphics[scale=0.14]{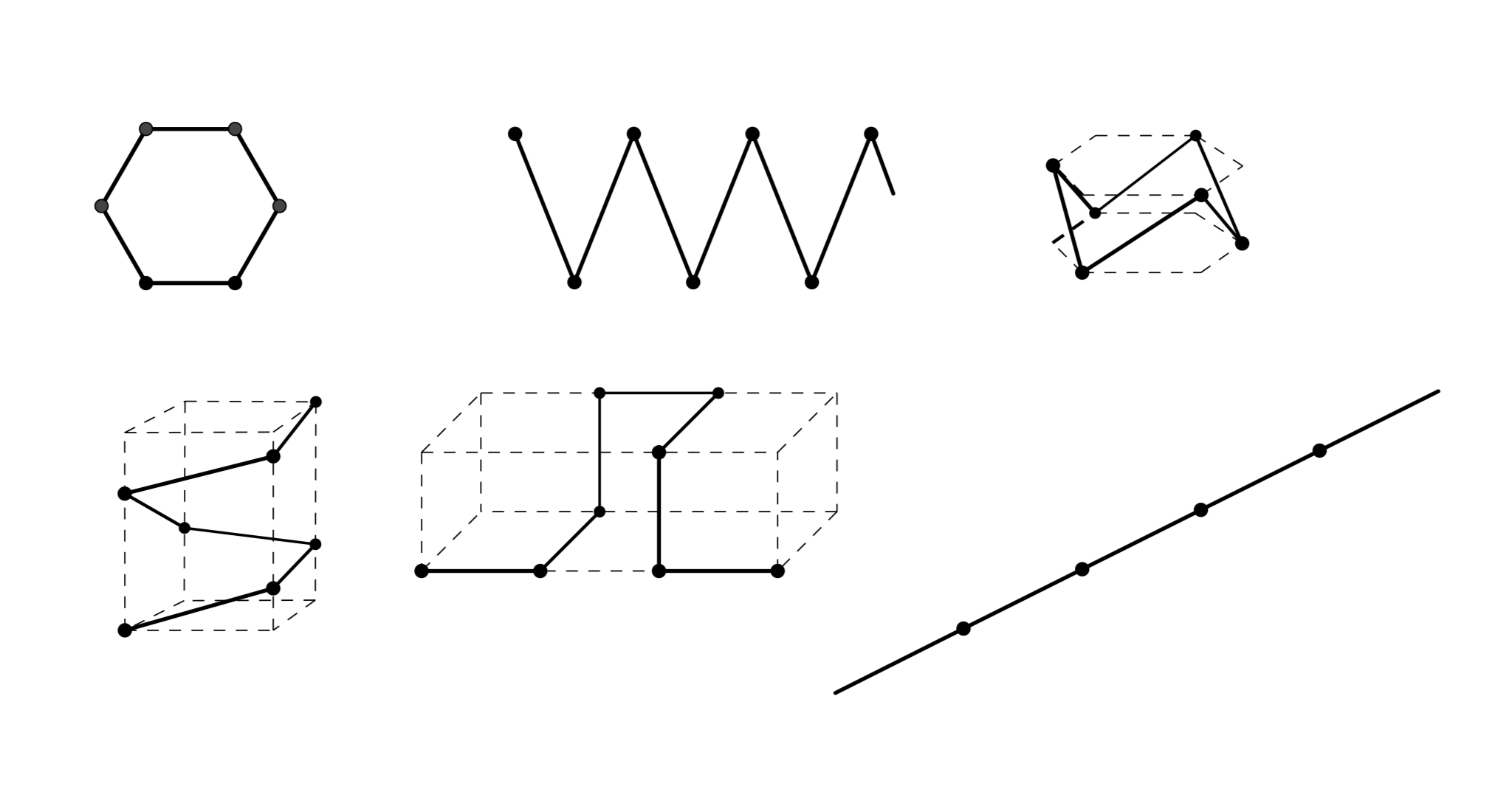}
\label{fig1}
\caption{Some polygons in $\mathbb{E}^3$.}
\end{figure}

The task of enumerating polyhedra in three-dimensional space according to its symmetry has proven to be difficult. Many efforts have been made over the last few years applying different techniques that have been diversifying. In 1977 B. Grünbaum \cite{grunbaum} introduces the following definition of polygon: A finite \textit{polygon} consists of a set of distinct points $\{v_1,\ldots ,v_n\}$ of $\mathbb{E}^d$ called \textit{the vertices} and a set of line segments $[v_iv_{i+1}]$ for $i=1,\ldots n-1$ and $[v_nv_1]$ in $\mathbb{E}^d$ called \textit{the edges}. Analogously, an infinite polygon can be defined in a similar way with an infinite set of points $\{\ldots,v_{-1},v_0,v_1,v_2,\ldots\}$ in such a way that each compact subset $K\subset \mathbb{E}^3$ must intersect a finite number of the edges of the polygon (Figure \ref{fig1}).\\

With this definition B. Grünbaum gives a classification of the regular polygons in $\mathbb{E}^3$ by means of the transitive action of the group of symmetries of the polygon. Each symmetry is a transformation of the space that preserves both the metric and the polygon itself. According to this, throughout this exposition we will work with the following definition: A geometric \textit{polyhedron} is a triplet $\mathcal{P}=(V,E,F)$ where $V$ is a set of points in $\mathbb{E}^3$ called \textit{the vertices}, $E$ is a set of segments called \textit{the edges} whose endpoints are elements of $V$, and $F$ is a set of polygons called \textit{the faces} whose edges are in $E$ and whose vertices are in $V$. Additionally $\mathcal{P}$ must satisfy the following properties:
\begin{enumerate}
\item Each edge of one of the faces is an edge of just one other face.
\item Each vertex of one face $f_i$ belongs to at least two other faces and all the faces which contain a vertex $v$ form a single 'circuit'; that is, they can be labelled cyclically so that neighbouring faces share an edge.
\item The family of polygons is connected: For any pair of edges $e,e'$ there exists a chain $e=e_0,f_1,e_1,f_2,\ldots , f_n,e_n=e'$ of edges and faces where each face $f_i$ contains $e_{i-1}$ and  $e_i$.
\item Each compact subset of $\mathbb{E}^3$ intersects $\mathcal{P}$ in a finite number of faces.
\end{enumerate}

Platonic solids are polyhedra having the property of being both convex and regular. The former is a geometrical property while the latter is related to the combinatorial structure of the polyhedron as detail below: A polyhedron has vertices, edges and faces. Let $v$ be a vertex, $e$ be an edge and $f$ be a face of a polyhedron $\mathcal{P}$. We say that $v$ \textit{is incident to} $e$, if $v$ is one of the endpoints of $e$. We say that $v$ \textit{is incident in} $f$ if $v$ is an endpoint of one of the edges in the polygon $f$. Finally we say that $e$ is incident to $f$ if $e$ is one of the edges in $f$. This can be summarized by saying that the incidence structure of a polyhedron $\mathcal{P}$ is a partial order in which the order relation is precisely the incidence, which we can denote by $\leq$. Any triplet $(v,e,f)$ such that $v\leq e\leq f$ is called \textit{a flag}. A symmetry of the polyhedron $\mathcal{P}$ is an isometry of the space that leaves $\mathcal{P}$ invariant. According to B. Grünbaum \cite{grunbaum}, a polyhedron $\mathcal{P}$ is said to be \textit{regular} if the group of symmetries of $\mathcal{P}$ acts transitively on the set of the flags. The regular polyhedra were classified by B. Grünbaum and A. W. Dress (see \cite{grunbaum} and \cite{dress}).

\section{Context of the problem.}

The purpose of this section is to document and put into context the problem of classifying all the fully transitive polyhedra. A polyhedron $\mathcal{P}=(V,E,F)$ is said to be \textit{fully transitive} if its group of symmetries acts transitively on all three sets: $V,E$ and $F$. Recently several efforts have been made towards the classification of abstract polytopes with a \textit{high degree of symmetry}, for example polyhedra whose group of symmetries has two orbits in the set of the flags, see for example \cite{hubard} where I. Hubard establishes seven different classes of 2-orbit abstract polytopes according to their incidence structure and the possible ways of organizing the flags in one or another orbit. She also establishes a classification of groups that can be automorphisms groups of abstract polytopes with two orbits. Four of these classes correspond to fully transitive polytopes. In \cite{schulte1} and \cite{schulte2}, E. Schulte classifies geometric polyhedra of one of these classes known as \textit{chiral polyhedra}. Regular polyhedra are, of course, those that have only one orbit in the set of the flags. There are also more general approaches about $k-$orbit polytopes. For these purposes it has been useful to abstract the combinatorial properties of geometric polyhedra giving rise to the theory of abstract polytopes. There are differences between this viewpoint and the geometric one. In this note we will focus on describing the geometric point of view of the following problem:

\begin{center}
\textit{Classify all geometric fully transitive polyhedra in $\mathbb{E}^3$.}
\end{center}

In 1988 Steven Lee Farris published his paper entitled \textit{Completely Classifying all vertex-transitive and edge-transitive polyhedra} \cite{farris_art} in which he establishes necessary conditions for fully transitive polyhedra and describes a method for generating them. These ideas are developed in his dissertation \textit{Fully transitive polyhedra} \cite{farris_tesis} under the supervision of B. Grünbaum, carrying out his method for the case of finite polyhedra. In 1990, John Merrick Eisenlohr obtained his PhD with his dissertation \textit{Fully-transitive polyhedra with crystallographic symmetry groups} \cite{eisenlohr} in which he applies Farris' ideas by means of an algorithm that generates all fully transitive polyhedra in any dimension $d\geq 2$. He also carries out the method for the case of planar infinite polyhedra. J. M. Eisenlohr establishes the terminology for the algorithm in which crystallographic groups play a central role and whose classification is a problem related to H. Poincaré's ideas about regular divisions of space, see for example \cite{goursat} and \cite{poincare}. He also mentions the topological context of the problem by studying the genus of the surface defined by the different polyhedra generated by the algorithm. In the next section we will describe the algorithm. 

\section{Brief description of the algorithm.}

The main problem is: Construct and classify all fully transitive polyhedra in $\mathbb{E}^3$. For this we will briefly discuss how crystallographic groups are defined. Let $E(3)$ be the group of isometries of $\mathbb{E}^3$. The set of 3-dimensional \textit{crystallographic} subgroups, denoted by $\mathcal{C}^3$ consists of the discrete subgroups of $E(3)$, i.e. subgroups that act discretely on $\mathbb{E}^3$. Equivalently we can say that a group $\Gamma \leq E(3)$ is a \textit{crystallographic group} if $\mathbb{E}^3/ \Gamma$ is compact (a good reference is \cite{ratcliffe}). Bieberbach's Theorem (see Theorem 7.5.3 in \cite{ratcliffe}) states that for each dimension there are only a finite number of isomorphism classes of crystallographic groups. It is well known that there are 17 different isomorphism classes of crystallographic groups in the plane. It is also known that there are 219 isomorphism classes in $\mathcal{C}^3$. A description of crystallography can be found in chapter 4 of the book \cite{coxeter}. The Eisenlohr - Farris algorithm starts by selecting a crystallographic group $\Gamma$ in order to generate all possible fully transitive geometric polyhedra with symmetry group $\Gamma$ and roughly consists of the following steps:
\begin{enumerate}
\item Characterize those vertex sets that are transitive, that is, that are obtained as the orbit of a point $p$ under the action of a crystallographic group.
\item For each vertex set $V$, list all graphs having vertices in $V$ that could serve as the 1-skeleton of a fully transitive polyhedron.
\item For each of these graphs, determine the different ways to fill in the faces to construct a fully transitive polyhedron.
\end{enumerate}

Throughout this exposition $p,u,v,w$ will represent points of $\mathbb{E}^3$, $\Gamma \in \mathcal{C}^3$ represents a crystallographic group of $\mathbb{E}^3$, $\Gamma (p)$ is the orbit of $p$ under the action of $\Gamma$. Let $V$ be the set $\Gamma (p)$ which is a discrete set of points of $\mathbb{E}^3$ on which $\Gamma$ acts transitively. We will sayy that $V$ is a \textit{crystallographic set} of points. Moreover every discrete set of points of $\mathbb{E}^3$ on which $\Gamma$ acts transitively can be obtained in this way (see \cite{eisenlohr}). Let $e$ be the line segment with endpoints $v_1,v_2\in V$. We construct a graph $E$ from $\Gamma (e)$, the orbit of $e$ under the action of $\Gamma$. In this way, $E$ is a fully transitive graph and furthermore every fully transitive graph can be constructed as the orbit of such an edge. Finally we will take a polygon $f$ with vertices and edges in $E$ and we will consider the orbit $\Gamma (f)$ under the action of $\Gamma$. In this way we obtain a fully transitive family $F$ of polygons that induces a fully transitive polyhedron \cite{eisenlohr}. Moreover, any fully transitive polyhedron can be constructed in this way (details can be found in \cite{eisenlohr}).

\section{Regular divisions of space and crystallographic sets of points.}

In \cite{poincare}, H. Poincaré studies the \textit{regular divisions of space} into an infinity of regions $R_0,R_1,\ldots ,R_i,\ldots$ such that each region $R_i$ can be obtained from the region $R_0$ by a transformation resulting from a composition of reflections. The classification of the crystallographic groups is related precisely to this idea: Each group $\Gamma \in \mathcal{C}^3$ consists of isometries. The \textit{Euclidean normalizer} of $\Gamma$, $N_E (\Gamma)$ is the subgroup $\{\kappa \in E(3) : \kappa \Gamma \kappa ^{-1}=\Gamma\}$. We can define a $\Gamma-$\textit{region} as the set $R_{\Gamma}=\mathbb{E}^3 / N_E (\Gamma)$. In this way the space can be covered with an infinite number of regions obtained from $R_{\Gamma}$ by a set of generators for $\Gamma$. Conversely, each regular division of the space defines a crystallographic group in $\mathcal{C}^3$.\\

J. M. Eisenlohr defines an equivalence relation by the euclidean similarity in $\mathcal{C}^3$ and also proves that for our purposes it is sufficient to consider one group for each similarity class. By defining $\mathcal{G}^3$ as the set containing one representative element of each similarity class of $\mathcal{C}^3$ we can reduce the number of groups to be considered. If $\mathcal{P}$ is a fully transitive polyhedron in $\mathbb{E}^3$, then there exists some group $\Gamma \in \mathcal{G}^3$, a \textit{base point} $u\in R_{\Gamma}$ and some polygon $f$ with vertices in $V=\Gamma (u)$ such that $\mathcal{P}$ is similar to the orbit $\Gamma (f)$. If we consider the set $\mathcal{O}=\{(\Gamma ,u):\Gamma \in \mathcal{G}^3, u\in R_{\Gamma}\}$ then a polyhedron is said to be generated by $(\Gamma, u)\in \mathcal{O}$ if it is the orbit of a polygon with vertices in $\Gamma (u)$. If $\mathcal{P} _1$ and $\mathcal{P} _2$ are polyhedra generated by $(\Gamma _1,u_1)$ and $(\Gamma _2,u_2)$ respectively and they happen to be similar polyhedra then one of the $\mathcal{P}_j$ has symmetry group larger than $\Gamma _j$. With this we are ruling out possible repetitions. The elements of $\mathcal{O}$ can be organized into \textit{cosymmetry classes}, a concept introduced by S. A. Robertson, S. Carter and H. R. Morton in \cite{robertson} that allows us to work by choosing one representative element $(\Gamma ,u)\in \mathcal{O}$ for each cosymmetry class since any two points in the same cosymmetry class will generate polyhedra with the same symmetry type (see \cite{eisenlohr}).
\\

Let us suppose for now that $\Gamma$ is the group defined by the \textit{honey-pie slice} $R_{\Gamma}\subset \mathbb{E}^3$ bounded by the planes $H_0,\ldots ,H_4$. The planes $H_0$ and $H_4$ are horizontal planes and the planes $H_1,H_2$ and $H_3$ are vertical planes with angles $\frac{\pi}{6},\frac{\pi}{2}$ and $\frac{\pi}{3}$ as shown in Figure 2. Clearly $R_{\Gamma}$ defines a regular division of the space and therefore the group $\Gamma$ generated by the reflections $\gamma _i$ through the planes $H_i$ is an element of $\mathcal{G}^3$. Our base point will be $u=H_0\cap H_1 \cap H_2$ as shown in the Figure 2.\\
If we restrict the action of $\Gamma$ to the plane $H_0$ we obtain a crystallographic set in the plane consisting of three \textit{lattice subsets}, which can be easily identified in the Figure 3. In general, the subgroup $\Lambda \subset \Gamma$ consisting of all translations in $\Gamma$ partitions the set $\Gamma (u)=V$ into \textit{lattice classes} $V=\bigcup _{i=1} ^k V_i$, where $v,v' \in V_i$ if and only if the translation by the vector $v-v'$ is an element of $\Gamma$.\\
By extending the action of the group to the whole space we will obtain several copies of this set in parallel planes generating a crystallographic set $V=\Gamma (u) \in \mathbb{E}^3$ with the same lattice classes.

\begin{figure}
\centering
\includegraphics[scale=0.14]{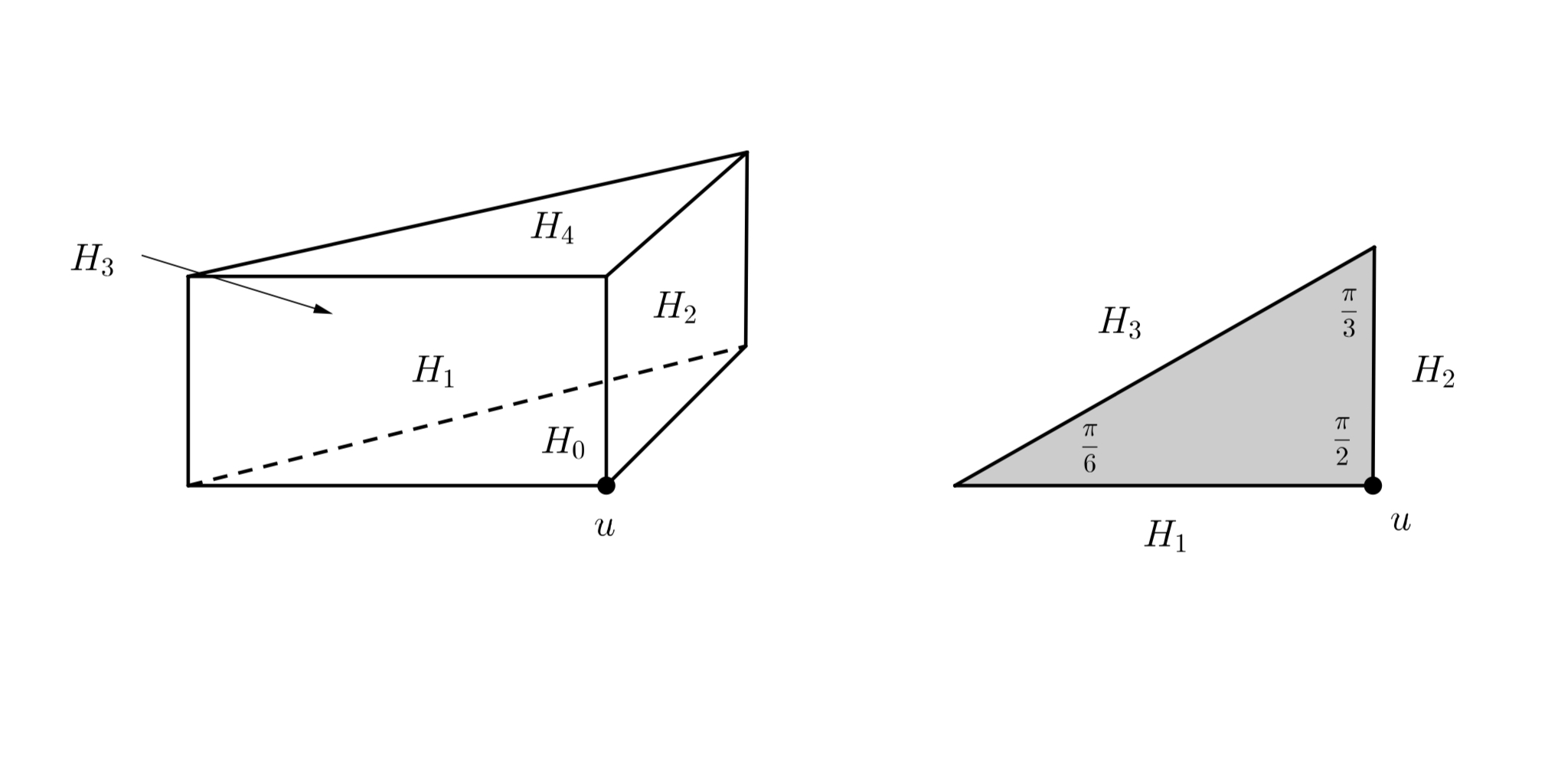}
\label{fig2}
\caption{A $\Gamma -$region.}
\end{figure}

\begin{figure}
\centering
\includegraphics[scale=0.14]{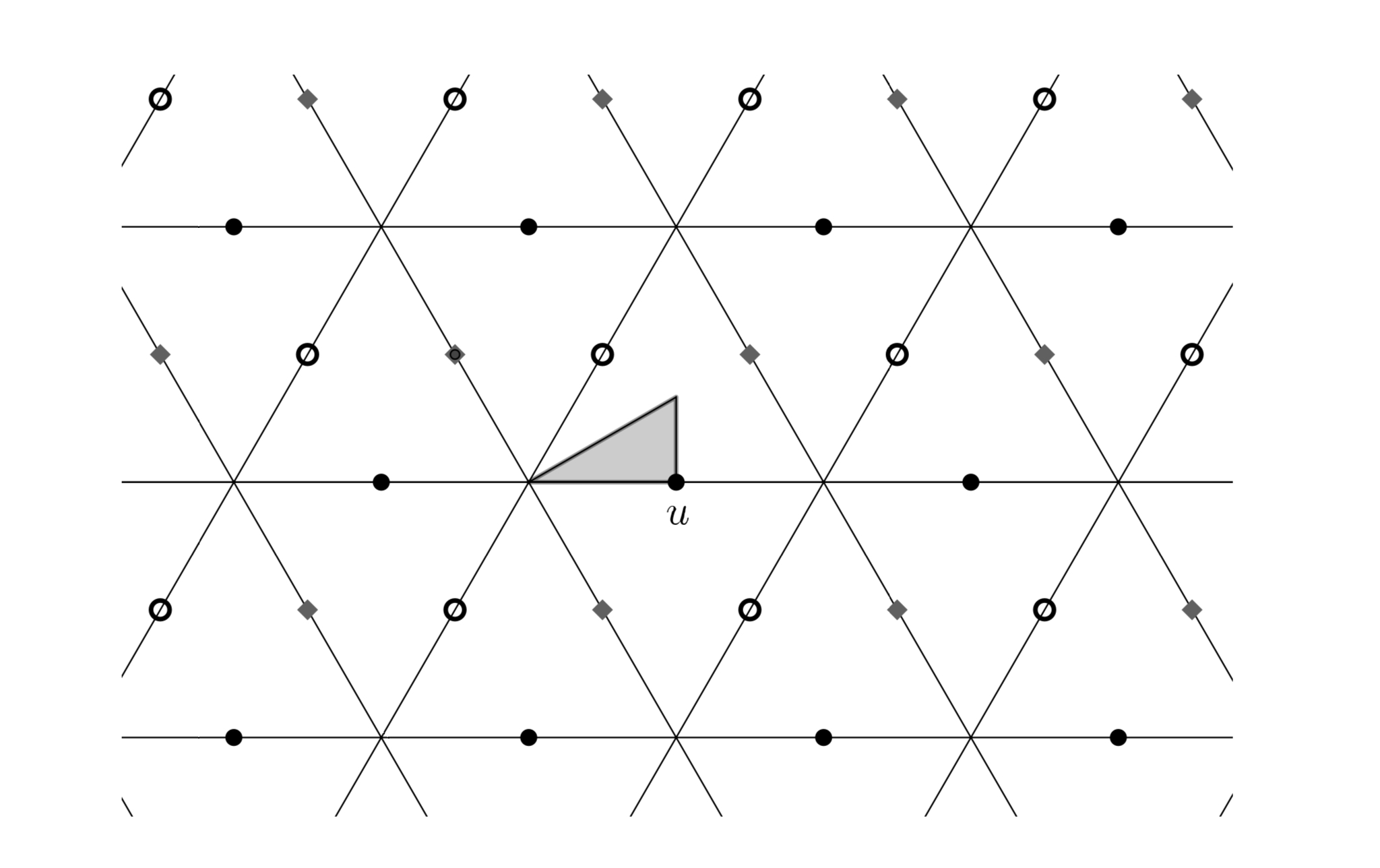}
\label{fig3}
\caption{Lattice subsets on the crystallographic set $V$.}
\end{figure}

\section{Edge sets and vertex figures.}

Now we will choose a \textit{base edge} with one of its endpoints being the base point $v$ and the other one being any other point $v$. Eisenlohr establishes that the structure of the set $Q_v$ of edges emanating from $u$ obtained from the edge $[uv]$ by letting act the stabilizer $\Gamma _u=\{\sigma \in \Gamma : \sigma (u)=u\}$ depends only on the lattice class in which the chosen vertex $v$ is located (see \cite{eisenlohr}, Proposition 1). We will call such a set $Q_v$ \textit{the star} of $v$. We will apply this idea to the set $V=\Gamma (u)$ of our previous example. In the plane $H_0$ we can identify the following points: $v=\gamma _2\gamma _3(u), w=\gamma _3(v),x=\gamma _1(w)$ and $y=\gamma _1(v)$. If we choose the edge $[uv]$ and by letting the stabilizer $\Gamma _u$ act on it, we obtain the following star $Q_v$ consisting of four edges emanating from $u$: $[uv],[uw],[ux]$ and $[uy]$ (Figure 4).\\

From this star we can now determine the possible vertex figures. The \textit{vertex figure} of a polyhedron $\mathcal{P}$ at a vertex $u$ of $\mathcal{P}$ is the polygon $[v_1\ldots v_q]$ where $v_1,\ldots ,v_q$ are the vertices of $\mathcal{P}$ adjacent to $u$, and each consecutive pair $v_k,v_{k+1}$, $1\leq k\leq q-1$, and $v_q, v_1$ belong to the same face. In this case we can produce three classes of vertex figure (Figure 5). \\

\begin{figure}
\centering
\includegraphics[scale=0.15]{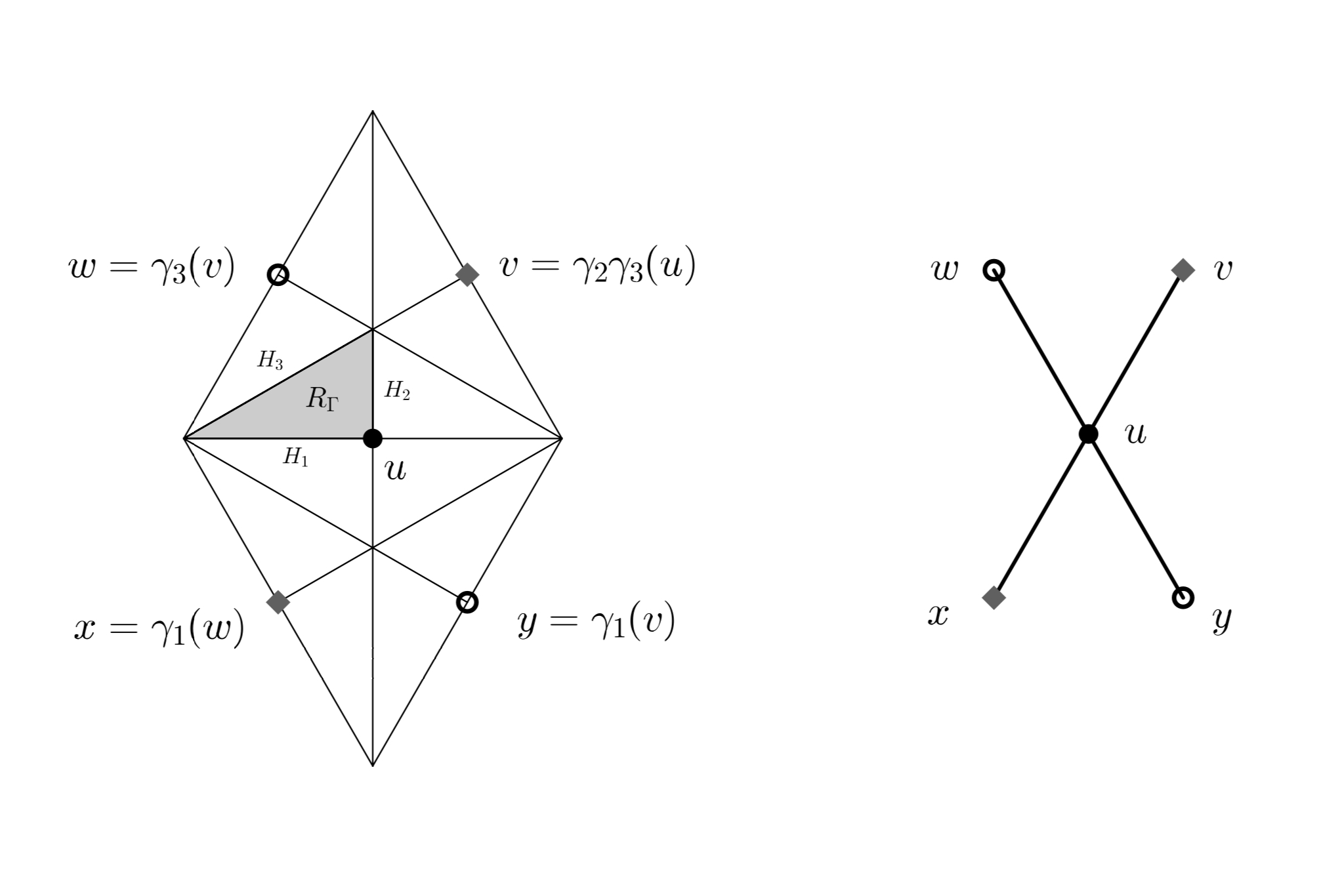}
\label{fig4}
\caption{The star of $v$ based at $u$.}
\end{figure}

\begin{figure}
\centering
\includegraphics[scale=0.15]{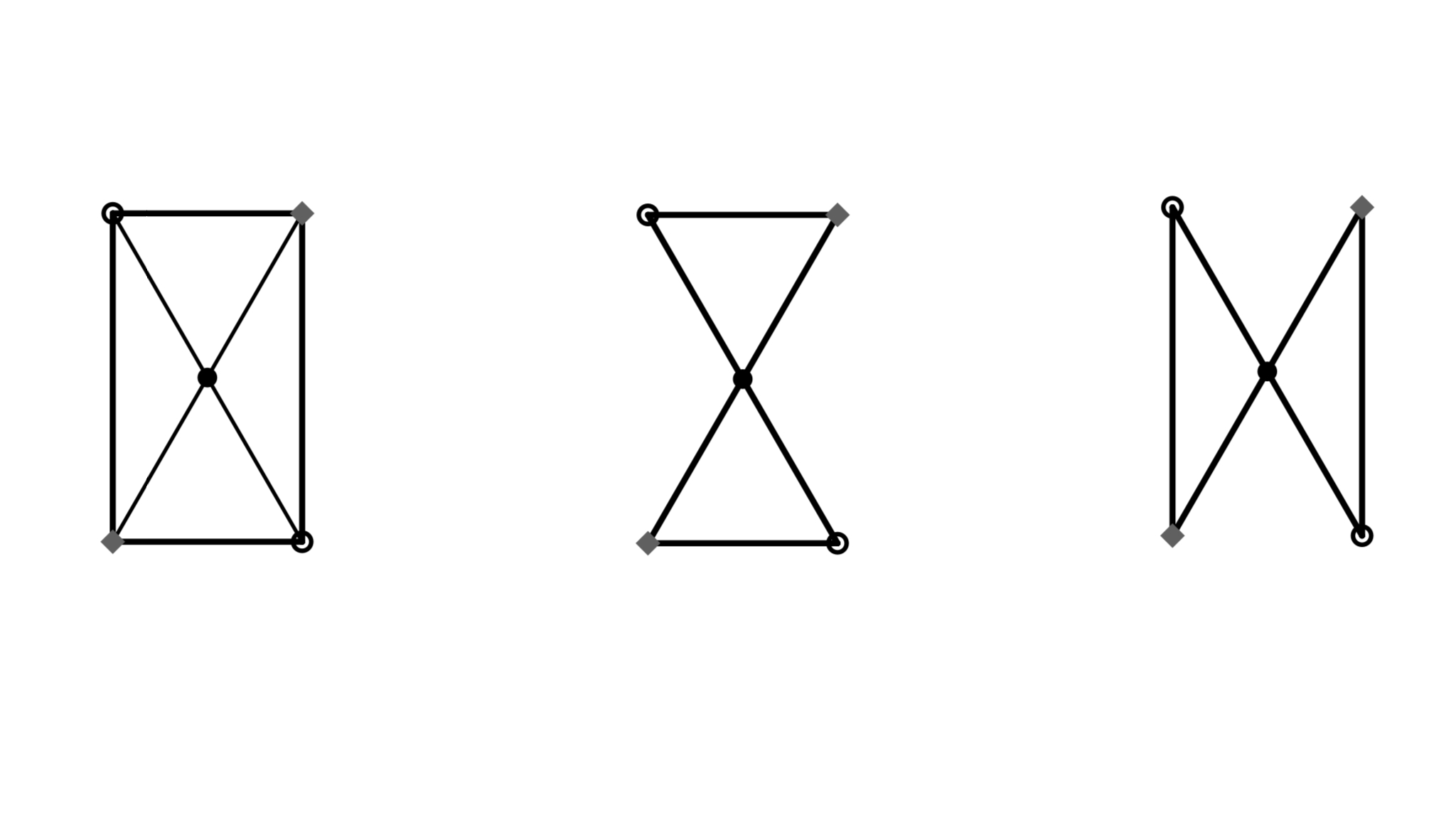}
\label{fig5}
\caption{Three possible vertex-figures.}
\end{figure}

S. L. Farris defines a \textit{face angle} of a polyhedron $\mathcal{P}$ as the planar angle $a$ between two consecutive edges of a face $f$ of $\mathcal{P}$, with $0\leq a \leq \pi$. Two angles $a_1$ and $a_2$ are called \textit{equivalent} if there is a symmetry of $\mathcal{P}$ which maps $a_1$ to $a_2$. The equivalence classes are called \textit{angle classes}. Farris' work establishes necessary conditions on the different angle classes of a fully transitive polyhedron:

\begin{theorem}\label{teo1} (Farris) Let $\mathcal{P}$ be any vertex-transitive and edge-transitive polyhedron, and let $v$ be any vertex of $\mathcal{P}$. Then one of the following statements is true:
\begin{enumerate}
\item $\mathcal{P}$ has exactly one angle class.
\item $\mathcal{P}$ has exactly two angle classes, and a circuit of angles at $v$ is $\alpha, \beta, \alpha, \beta ,\ldots ,\beta$.
\item $\mathcal{P}$ has exactly three angle classes, and a circuit of angles at $v$ is $\alpha, \beta _1, \alpha, \beta _2 ,\ldots ,\beta _2$
\end{enumerate}
\end{theorem}

\begin{theorem}\label{teo2} (Farris) Let $\mathcal{P}$ be a fully-transitive polyhedron, and let $f$ be any face of $\mathcal{P}$. One of the following is true:
\begin{enumerate}
\item $\mathcal{P}$ has exactly one angle class.
\item $\mathcal{P}$ has exactly two angle classes, and a circuit of angles of $f$ is $\alpha, \beta, \alpha, \beta ,\ldots ,\beta$ if $f$ is finite, and $\alpha, \beta, \alpha, \beta ,\ldots $ if $f$ is infinite.
\item $\mathcal{P}$ has exactly three angle classes, and a circuit of angles of $f$ is $\alpha, \beta _1, \alpha, \beta _2 ,\ldots ,\beta _2$ if $f$ is finite, and $\alpha, \beta _1, \alpha, \beta _2 ,\ldots $ if $f$ is infinite
\end{enumerate}
\end{theorem}

\begin{lemma}\label{lem3} (Eisenlohr) Suppose we have fixed a vertex set $V$, an edge set $E$ and a vertex figure $\Sigma$. Then if there are two or three angle classes in $\Sigma$, then there is at most one fully-transitive polyhedron with vertex set $V$, edge set $E$ and vertex figure $\Sigma$.
\end{lemma}

With this information (see \cite{farris_tesis} and \cite{eisenlohr}) we can determine the possible ways to fill in the faces for each of these vertex figures. This is done by J. M. Eisenlohr in his dissertation, for in this case we have obtained plane vertex figures and therefore we will obtain one of the fully transitive polyhedra in the plane. In the following section we give an example that illustrates how the algorithm can be applied in three-dimensional space.

\section{An example in three-dimensional space.}

Let us consider the planes $\gamma _4(H_0)$ and $\gamma _0\gamma _4(H_0)$ parallel to $H_0$. In each of these planes is located a copy of the crystallographic set we have analysed in the preceding paragraphs, so that for each vertex in that set there is a copy at the corresponding level. We will mark with sub index 1 the corresponding images of the points of $H_0 \cap V$ that are in the plane $\gamma _4(H_0)$ and with 2 those in the plane $\gamma _0\gamma _4(H_0)$. Let the stabilizer $\Gamma _u$ act on the edge $uv_1$ in order to obtain the star $Q_{v_1}$ as shown in Figure 6.\\

\begin{figure}
\centering
\includegraphics[scale=0.14]{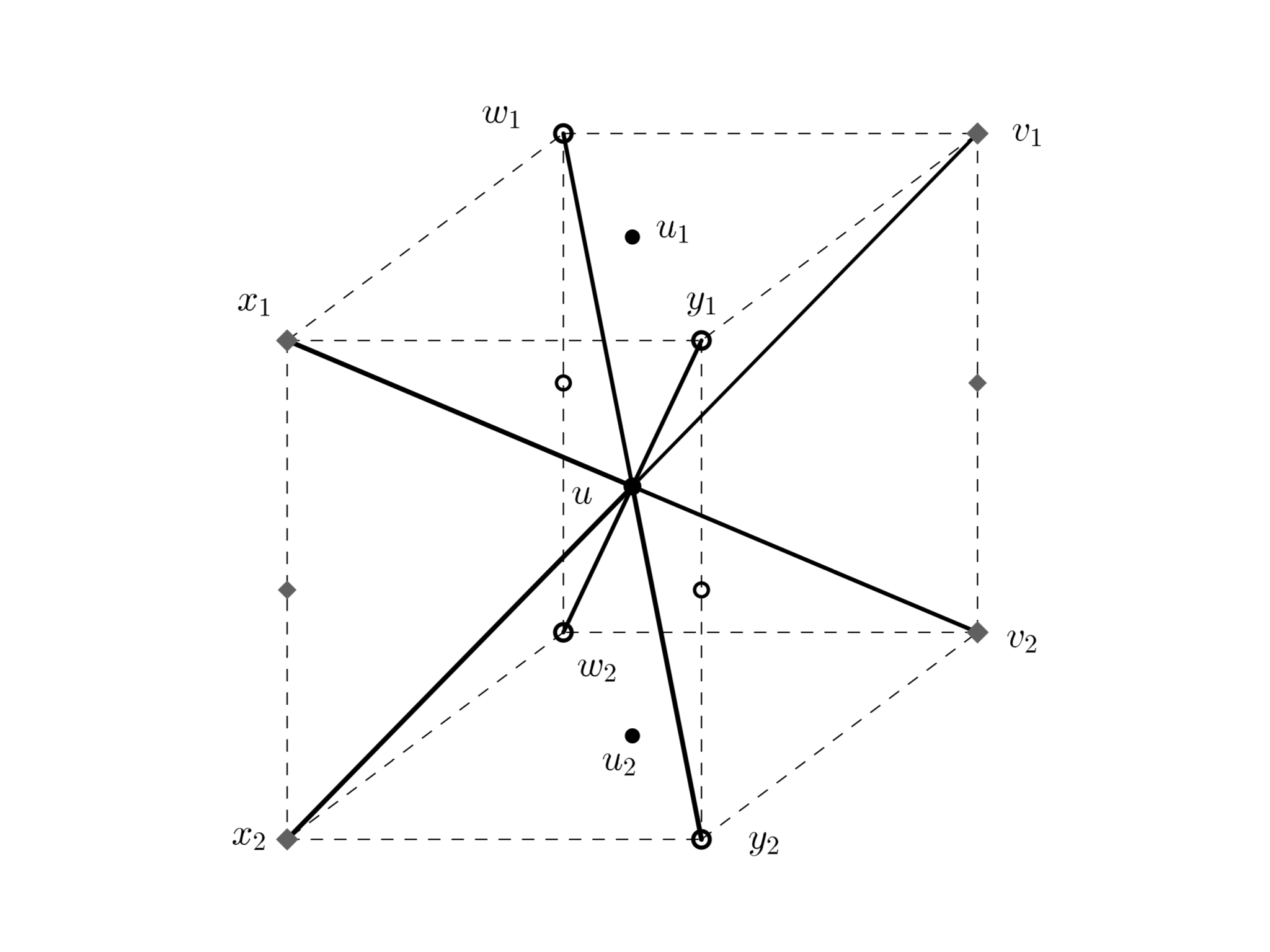}
\label{fig6}
\caption{The star of $v_1$ based at $u$.}
\end{figure}

\begin{figure}
\centering
\includegraphics[scale=0.14]{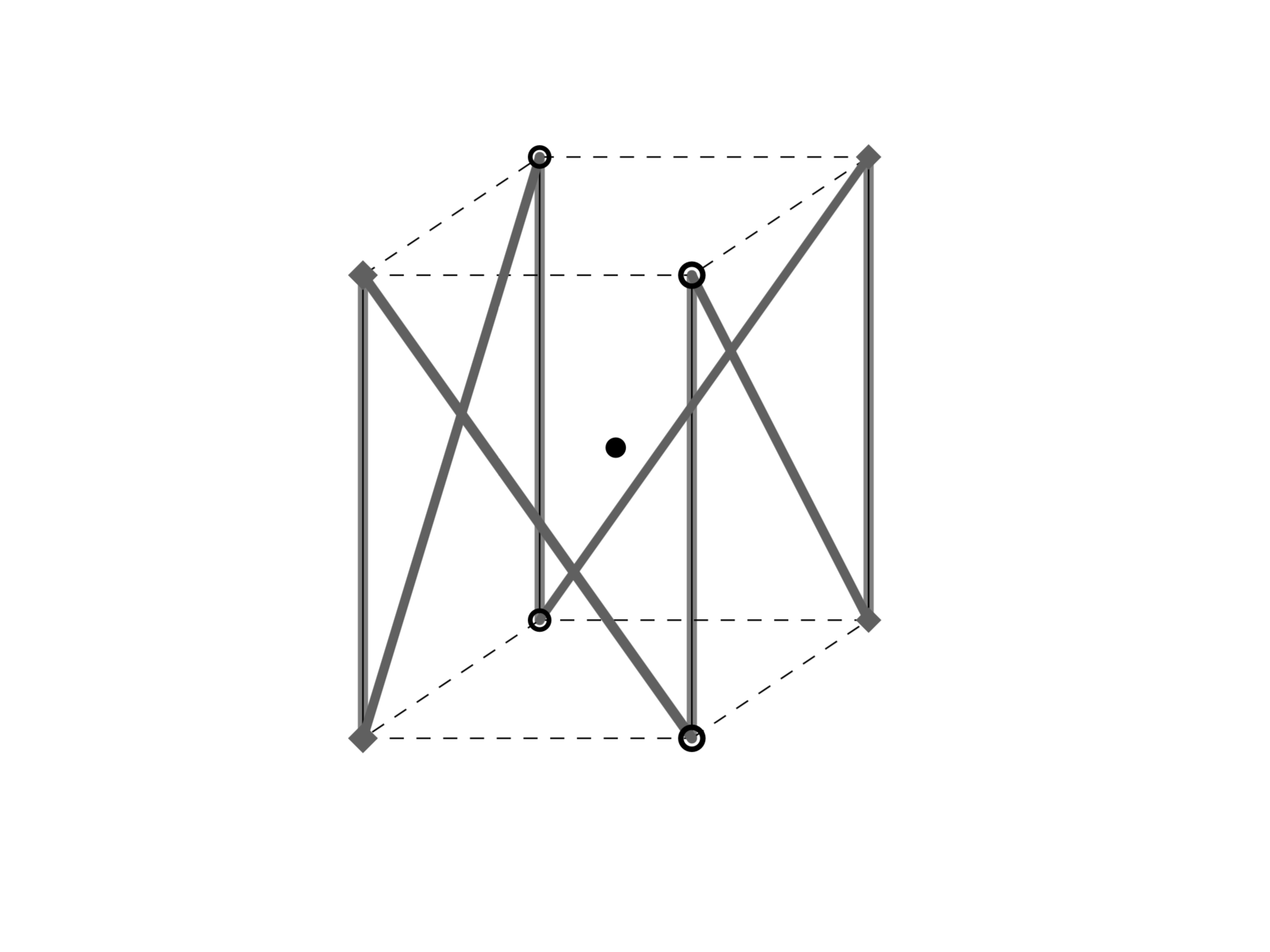}
\label{fig7}
\caption{An example of vertex-figure.}
\end{figure}

Now, by using Theorem \ref{teo1} and Theorem \ref{teo2} we can determine the different vertex figures for this star. Let's consider for example the Figure 7.\\
In this case we have three angle classes defined by $\alpha =\angle x_1 ux_2$, $\beta _1=\angle x_1uy_2$ and $\beta _2 =\angle y_1uv_2$ so, by Lemma \ref{lem3} there is only one way to fill in the faces. In order to define a face $f$ we use the proof of Lemma \ref{lem3} (Lemma 2.9 in \cite{eisenlohr}): Since the vertex figure is specified whenever we have chosen the first 2 edges defining certain angle, there are two choices for the next edge, but by means of Theorem \ref{teo2} (Theorem 2.5 of \cite{farris_art} we must choose the next edge so that the face angles alternate between the alpha and the beta angles, and since the angles alternate in the vertex figure, by applying Theorem \ref{teo1} (Theorem 2.4 of \cite{farris_art}) there is only one choice. With these ideas in mind it is not too difficult to see that $f$ is a \textit{zig-zag spiral} described as follows: Let $x$ be any vertex of the crystallographic set $V$ located at, say, $H_0$. Let's imagine we describe a downward spiral from $x$. The first direction will be $[xy]$ in back and forth way: $[xy]$ and $[yx]$ descending two levels. Then the next direction will be ]$[xu]$ and $[ux]$ descending two levels. Then $[xz]$ and $[zx],[xw].[wx]$ and again $[xy],[yx]$ and so on (Figure 8).
 
 \begin{figure}
\centering
\includegraphics[scale=0.17]{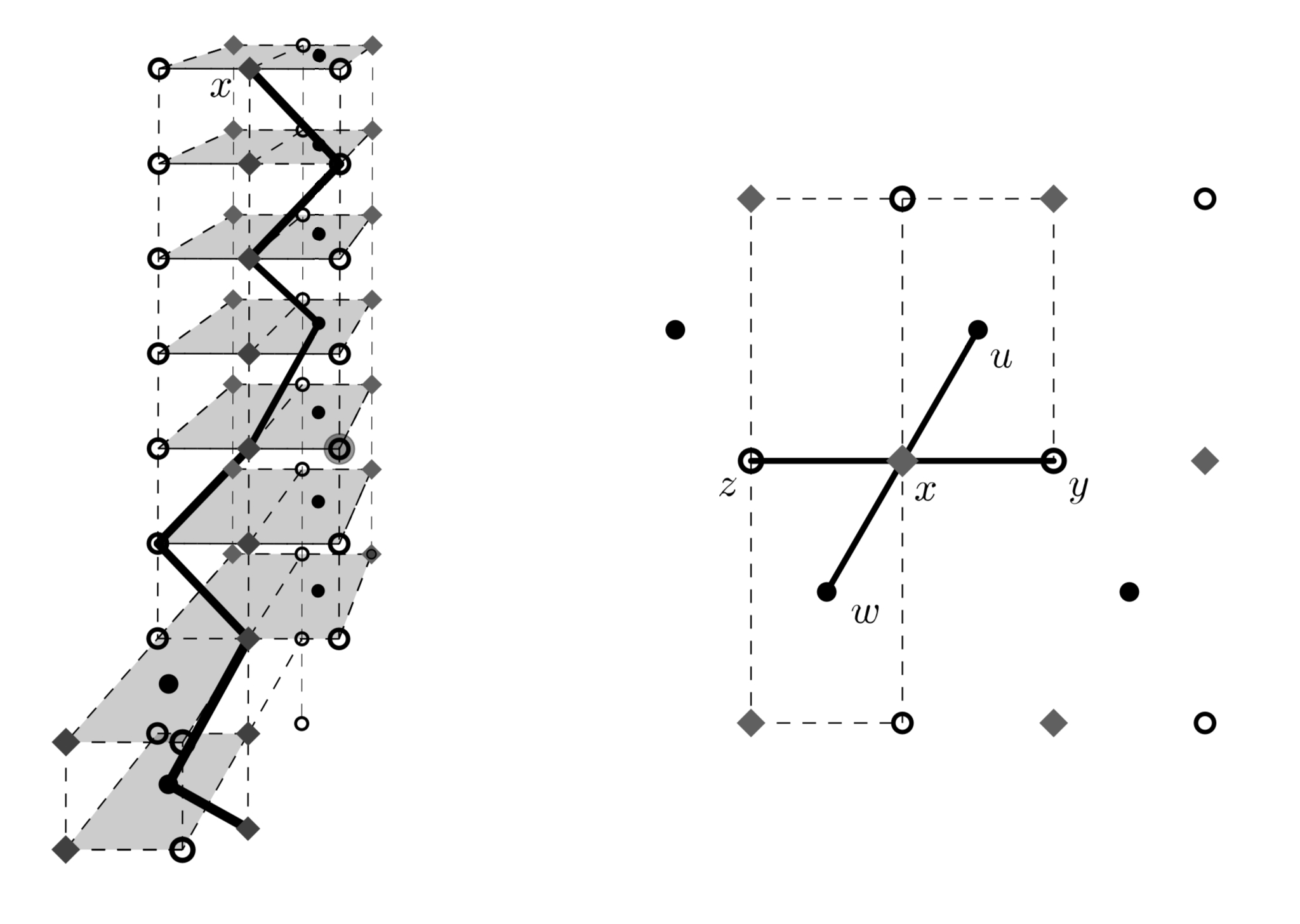}
\label{fig8}
\caption{Zig-zag spiral as a face.}
\end{figure}

In view of Eisenlohr and Farris results this polygon defines an hexagonal \textit{zig-zag spiralhedron} $\mathcal{S}_1$, a polyhedron whose group of symmetries is hexagonal, according to Eisenlohr's terminology. Furthermore:

\begin{prop} This hexagonal zig-zag spiralhedron $S_1$ is a fully transitive polyhedron.
\end{prop}

\vspace{1cm}

\noindent \textit{Funding and competing interests: Project supported by postdoctoral fellowship from DGAPA-UNAM. The datasets generated during the current study are available from the corresponding
author on reasonable request.}

\end{document}